\documentclass[number,dvips]{arxbj}
\usepackage{mathbh}

\aid{0}
\volume{13}
\issue{4}
\pubyear{2007}
\firstpage{1000}
\lastpage{1022}
\doi{10.3150/07-BEJ6044}

\makeatletter

\newtheorem{theo}{Theorem}
\newtheorem{coro}{Corollary}
\newtheorem{lem}{Lemma}
\newproclaim{ma}{(MA)}
\newproclaim{mah}{(MAH)}
\newtheorem{pro}{Proposition}
\newremark{rem}{Remark}
\newproclaim{defi}{Definition}
\newcommand{\1}{\mathbh{1}}
\makeatother

\begin{document}
\begin{frontmatter}

\title{Optimal rates of aggregation in classification under low noise
assumption}
\runtitle{Aggregation of classifiers}

\begin{aug}
\author[a]{\fnms{Guillaume} \snm{Lecu\'e}\corref{}\ead[label=e1]{lecue@ccr.jussieu.fr}}
\runauthor{G. Lecu\'e}
\pdfauthor{Guillaume Lecue}
\address[a]{Laboratoire de Probabilit\'es et Mod\`eles
Al\'eatoires (UMR CNRS 7599), Universit\'e Paris
VI,\\ 4 pl. Jussieu, BP 188,
75252 Paris, France. \printead{e1}}
\end{aug}

\received{\smonth{3} \syear{2006}}
\revised{\smonth{4} \syear{2007}}

%
\begin{abstract}
In the same spirit as Tsybakov, we define the optimality of
an aggregation procedure in the problem of classification. Using
an aggregate with exponential weights, we obtain an optimal rate
of convex aggregation for the hinge risk under the margin
assumption. Moreover, we obtain an optimal rate of model selection
aggregation under the margin assumption for the excess Bayes risk.
\end{abstract}

%

\end{frontmatter}
%
\section{Introduction}
Let $(\mathcal{X},\mathcal{A})$ be a measurable space. We consider a random
variable $(X,Y)$ on $\mathcal{X}\times\{-1,1\}$ with probability
distribution denoted by $\pi$. Denote by $P^X$ the marginal of
$\pi$ on $\mathcal{X}$ and by $\eta(x)\stackrel{\mathrm
{def}}{=}\mathbb{P}(Y=1|X=x)$ the
conditional probability function of $Y=1$, knowing that $X=x$. We
have $n$ i.i.d. observations of the couple $(X,Y)$ denoted by
$D_n=((X_i,Y_i))_{i=1,\ldots,n}$. The aim is to predict the output
label $Y$ for any input $X$ in $\mathcal{X}$ from the observations $D_n$.

We recall some usual notation  for the classification
framework. A \textbf{prediction rule} is a measurable function
$f\dvtx\mathcal{X}\longmapsto\{-1,1\}$. The \textbf{misclassification error}
associated with $f$ is
\[
R(f) = \mathbb{P}\bigl(Y\neq f(X)\bigr).
\]
It is well
known (see, e.g., Devroye \textit{et al.} \cite{dgl96})
that
\[
\min_{f:\mathcal{X}\longmapsto\{-1,1\}} R(f) = R(f^*)\stackrel
{\mathrm{def}}{=}
R^{*},
\]
where the prediction rule $f^*$, called the \textbf{Bayes rule}, is defined by
\[
f^*(x)\stackrel{\mathrm{def}}{=}{\rm{sign}}\bigl(2\eta(x)-1\bigr)\qquad \forall
x\in\mathcal{X}.
\]
The minimal risk $R^*$ is called the \textbf{Bayes risk}. A
\textbf{classifier} is a function, $\hat{f}_n=\hat{f}_n(X,D_n)$,
measurable with respect to $D_n$ and $X$ with values in $\{-1,1\}$,
that assigns to the sample $D_n$ a prediction rule
$\hat{f}_n(\cdot,D_n)\dvtx\mathcal{X}\longmapsto\{-1,1\}$. A key
characteristic of
$\hat{f}_n$ is the \textbf{ generalization error} $\mathbb{E}[R(\hat{f}_n)]$,
where
\[
R(\hat{f}_n) \stackrel{\mathrm{def}}{=}\mathbb{P}\bigl(Y\neq\hat
{f}_n(X) |D_n\bigr).
\]
The aim of
statistical learning is to construct a classifier $\hat{f}_n$ such that
$\mathbb{E}[R(\hat{f}_n)]$ is as close to $R^*$ as possible. Accuracy
of a
classifier $\hat{f}_n$ is measured by the value $\mathbb{E}[R(\hat
{f}_n) - R^*]$,
called the \textbf{excess Bayes risk} of $\hat{f}_n$. We say that the classifier
$\hat{f}_n$ learns with the convergence rate $\psi(n)$, where
$(\psi(n))_{n\in\mathbb{N}}$ is a decreasing sequence, if there
exists an absolute constant $C>0$ such that for any integer $n$,
$\mathbb{E}[R(\hat{f}_n) - R^*]\leq C\psi(n)$.

Given a convergence rate, Theorem 7.2 of Devroye  \textit{et al.} \cite{dgl96} shows that no
classifier can learn at least as fast as this rate for any
arbitrary underlying probability distribution $\pi$. To achieve
rates of convergence, we need a complexity assumption on the set
which the Bayes rule $f^*$ belongs to. For instance,
Yang \cite{yang99,yang299} provide examples of classifiers learning
with a given convergence rate under complexity assumptions. These
rates cannot be faster than $n^{-1/2}$ (cf. Devroye  \textit{et al.} \cite{dgl96}).
Nevertheless, they can be as fast as $n^{-1}$ if we add a control on
the behavior of the conditional probability function $\eta$ at the
level $1/2$ (the distance $|\eta(\cdot)-1/2|$ is sometimes called
the \textbf{margin}). For the problem of
discriminant analysis, which is close to our classification problem, Mammen  and Tsybakov \cite{mt99}
and Tsybakov  \cite{tsy04} have introduced the following assumption.

\begin{ma*}[Margin (or low noise) assumption] The
probability distribution $\pi$ on the space $\mathcal{X}\times\{-1,1
\}$
satisfies $\operatorname{MA}(\kappa $) with $1\leq\kappa <+\infty$ if there exists $c_0>0$
such that
%
\begin{equation}\label{e4}
\mathbb{E}[ |f(X)-f^{*}(X)|]\leq c_0 \bigl(R(f)-R^{*}
\bigr)^{1/\kappa },
\end{equation}
for any measurable function $f$ with values in $\{-1,1\}.$
\end{ma*}

 According to Tsybakov \cite{tsy04} and Boucheron \textit{et al.} \cite{bbl05}, this
assumption is equivalent to a control on the margin given by
\[
\mathbb{P}[|2\eta(X)-1|\leq t]\leq c t^{\alpha}\qquad\forall0\leq t <1.
\]
Several example of \textbf{fast rates}, that is, rates faster than
$n^{-1/2}$, can be found in Blanchard \textit{et al.}
\cite{bbm04}, Steinwart and Scovel \cite{ss05,ss04},
Massart \cite{m00}, Massart and N\'{e}d\'{e}lec \cite{mn03},
Massart \cite{m04} and Audibert and Tsybakov \cite{at05}.

%

The paper is organized as follows. In Section, \ref{sectiondef} we
introduce definitions and   procedures which are used throughout
the paper. Section \ref{sectionoptimalhinge} contains oracle
inequalities for our aggregation procedures w.r.t. the excess hinge
risk. Section \ref{sectionExcessRisk} contains similar results for
the excess Bayes risk. Proofs are postponed to Section \ref{proofs}.

\section{Definitions and procedures}\label{sectiondef}
\subsection{Loss functions}
Convex surrogates $\phi$ for the classification loss are often used
in algorithm (Cortes and Vapnic \cite{cv95}, Freund and Schapire \cite{fs97},
Lugosi and Vayatis \cite{lv04}, Friedman \textit{et al.} \cite{fht00},
B\"{u}hlman and Yu \cite{by02}, Bartlett \textit{et al.} \cite{bfls98,bjm03}).
Let us introduce some notation. Take $\phi$ to be a measurable function
from $\mathbb{R}$ to $\mathbb{R}$. The risk associated with the loss
function $\phi$ is called the ${\mathbf{\phi}}$-{\bf risk} and is
defined by
\[
A^{(\phi)}(f)\stackrel{\mathrm{def}}{=}\mathbb{E}[\phi(Yf(X))],
\]
where
$f\dvtx\mathcal{X}\longmapsto\mathbb{R}$ is a measurable function. The
\textbf{empirical} $\phi$-\textbf{risk} is defined by
\[
A^{(\phi)}_n(f)\stackrel{\mathrm{def}}{=}\frac{1}{n}\sum
_{i=1}^n\phi(Y_if(X_i))
\]
and we
denote by $A^{(\phi)*}$ the infimum over all real-valued functions
$\inf_{f:\mathcal{X}\longmapsto\mathbb{R}}{A^{(\phi)}(f)}$.

Classifiers obtained by minimization of the empirical $\phi$-risk,
for different convex losses, have been proven to have very good
statistical properties (cf. Lugosi and Vayatis \cite{lv04}, Blanchard \textit{et al.} \cite{blv03},
Zhang \cite{z04}, Steinwart and Scovel \cite{ss05,ss04}
and Bartlett \textit{et al.} \cite{bjm03}). A wide variety of classification methods in
machine learning are based on this idea, in particular, on using
the convex loss $\phi(x)\stackrel{\mathrm{def}}{=}\max(1-x,0)$
associated with
support vector machines (Cortes and Vapnik~\cite{cv95}, Sch\"olkopf and Smola \cite{ssbook02}), called the
\textbf{hinge loss}. The corresponding risk is called the \textbf{hinge
risk} and is defined by
\[
A(f)\stackrel{\mathrm{def}}{=}\mathbb{E}\bigl[\max\bigl(1-Yf(X),0\bigr)\bigr],
\]
for any measurable function
$f\dvtx\mathcal{X}\longmapsto\mathbb{R}$. The \textbf{optimal hinge
risk} is
defined by
%
\begin{equation}\label{OHR}
A^*\stackrel{\mathrm{def}}{=}\inf_{f:\mathcal{X}\longmapsto\mathbb
{R}} A(f).
\end{equation}
It is easy to check that the Bayes rule $f^*$ attains the infimum in
(\ref{OHR}) and that
%
\begin{equation}\label{e1}
R(f)-R^*\leq A(f)-A^*,
\end{equation}
for any measurable function $f$ with values in $\mathbb{R}$ (cf.
Lin~\cite{l99} and generalizations in Zhang \cite{z04} and Bartlett \textit{et al.} \cite{bjm03}),
where we extend the definition of $R$ to the class of real-valued
functions by $R(f)=R(\mathrm{sign}(f))$. Thus, minimization of the
\textbf{excess hinge risk}, $A(f)-A^*$, provides a reasonable
alternative for minimization of the excess Bayes risk, $R(f)-R^*$.

\subsection{Aggregation procedures}\label{subsectionAggregationProcedures}
Now, we introduce the problem of aggregation and the aggregation
procedures which will be studied in this paper.

Suppose that we have $M\geq2$ different classifiers
$\hat{f}_1,\ldots,\hat{f}_M$ taking values in $\{-1,1\}$. The
problem of model selection type aggregation, as studied in
Nemirovski \cite{n00}, Yang \cite{yang00}, Catoni \cite{cat99,catbook01} and Tsybakov \cite{tsy03}, consists
of the
construction of a new classifier $\tilde{f}_n$ (called an
\textbf{aggregate}) which approximately mimics the best classifier
among $\hat{f}_1,\ldots,\hat{f}_M$. In most of these papers the
aggregation is based on splitting  the sample into two independent
subsamples, $D^{1}_{m}$ and $D^{2}_{l}$, of sizes $m$ and $l$,
respectively, where $m+l=n$. The first subsample, $D_{m}^{1}$, is used
to construct the classifiers $\hat{f}_1,\ldots,\hat{f}_M$ and the
second subsample, $D_{l}^{2}$, is used to aggregate them, that is to
construct a new classifier that mimics, in a certain sense, the
behavior of the best among the classifiers $\hat{f}_j,j=1,\ldots,M$.

In this paper, we will not consider the sample splitting and will
concentrate only on the construction of aggregates (following
Juditsky and Nemirovski \cite{jn00}, Tsybakov \cite{tsy03}, Birg{\'e} \cite{b04},
Bunea \textit{et al.} \cite{btw105}). Thus, the first subsample is
fixed and, instead of classifiers $\hat{f}_1,\ldots,\hat{f}_M$, we
have fixed prediction rules $f_{1},\ldots,f_M $. Rather than working
with a part of the initial sample we will suppose, for notational
simplicity, that the whole sample $D_n$ of size $n$ is used for the
aggregation step instead of a subsample $D_l^2$.

Let $\mathcal{F}=\{f_1,\ldots,f_M\}$ be a finite set of real-valued
functions, where $M\geq2$.
An \textbf{aggregate} is a real-valued statistic of the form
\[
\tilde{f}_n=\sum_{f\in\mathcal{F}}w^{(n)}(f)f,
\]
where the weights
$(w^{(n)}(f))_{f\in\mathcal{F}}$ satisfy
\[
w^{(n)}(f)\geq0 \quad\mbox{and}\quad\sum_{f\in\mathcal{F}}w^{(n)}(f)=1.
\]
Let $\phi$ be a convex loss for classification.
The Empirical Risk Minimization aggregate (\textbf{ERM}) is defined by the weights
\[
w^{(n)}(f)=
\cases{
\displaystyle1, &\quad  for one  $f\in\mathcal{F}$  such that
$A_n^{(\phi)}(f)=\min\limits_{g\in\mathcal{F}}A_n^{(\phi)}(g)$,\cr
0, &\quad  for all other $f\in\mathcal{F}$,
}
\qquad\forall f\in\mathcal{F}.
\]
The ERM aggregate is denoted by $\tilde{f}_n^{(\mathrm{ERM})}$.

The \textbf{averaged ERM} aggregate is defined by the weights
\[
w^{(n)}(f)=\cases{
1/N, &\quad  if  $A_n^{(\phi)}(f)=\min\limits_{g\in\mathcal{F}}A_n^{(\phi
)}(g)$,\cr
0, &\quad otherwise,
}
\qquad\forall f\in\mathcal{F},
\]
where $N$ is the number of functions in $\mathcal{F}$ minimizing the
empirical $\phi$-risk.
The averaged ERM aggregate is denoted by $\tilde{f}_n^{(\mathrm{AERM})}$.

The Aggregation with
Exponential Weights aggregate (\textbf{AEW}) is defined by the weights
%
\begin{equation}\label{AEW}w^{(n)}(f)=\frac{\exp( -nA_n^{(\phi
)}(f))}{\sum_{g\in\mathcal{F}}\exp(
-nA_n^{(\phi)}(g))}\qquad\forall f\in\mathcal{F}.
\end{equation}
The AEW aggregate is denoted by $\tilde{f}_n^{(\mathrm{AEW})}$.

The \textbf{cumulative AEW} aggregate is an on-line procedure defined by
the weights
\[
w^{(n)}(f)=\frac{1}{n}\sum_{k=1}^{n}\frac{\exp( -kA_k^{(\phi
)}(f))}{\sum_{g\in\mathcal{F}}\exp(
-kA_k^{(\phi)}(g))}\qquad\forall f\in\mathcal{F}.
\]
The cumulative
AEW aggregate is denoted by $\tilde{f}_n^{(\mathrm{CAEW})}$.

When $\mathcal{F}$ is a class of prediction rules, intuitively, the AEW
aggregate is more robust than the ERM aggregate w.r.t. the problem
of overfitting. If the classifier with smallest empirical risk is
overfitted, that is, if it fits too many to the observations, then the ERM
aggregate will be overfitted. But, if other classifiers in $\mathcal
{F}$ are
good classifiers, then the aggregate with exponential weights will
consider their ``opinions'' in the final decision procedure and these
opinions can balance with the opinion of the overfitted classifier
in $\mathcal{F}$, which can be false because of its overfitting
property. The
ERM only considers the ``opinion'' of the classifier with the smallest
risk, whereas the AEW takes into account all of the opinions of the
classifiers in the set $\mathcal{F}$.

The exponential weights, defined in (\ref{AEW}), can be found in
several situations. First, one can check that the solution of the
minimization problem
%
\begin{equation}\label{equationminimization}\min\Biggl(\sum_{j=1}^M
\lambda_jA_n^{(\phi)}(f_j)+\epsilon\sum_{j=1}^M\lambda_j\log
\lambda_j\dvtx\sum_{j=1}^M\lambda_j\leq1,
\lambda_j\geq0,j=1,\ldots,M\Biggr)
\end{equation}
for all $\epsilon>0$
is
\[
\lambda_j=\frac{\exp(-(A_n^{(\phi)}(f_j))/{\epsilon}
)}{\sum_{k=1}^M
\exp(-(A_n^{(\phi)}(f_k))/{\epsilon})}\qquad \forall
j=1,\ldots,M.
\]
Thus, for $\epsilon=1/n$, we find the exponential
weights used for the AEW aggregate. Second, these weights can also
be found in the theory of prediction of individual sequences (cf. Vovk
\cite{v90}).

\subsection{Optimal rates of aggregation}
Now, we introduce a concept of optimality for an aggregation
procedure and for rates of aggregation, in the same spirit as in
Tsybakov \cite{tsy03} (where the regression problem is treated). Our aim is
to prove that the aggregates introduced above are optimal in the
following sense. We denote by $\mathcal{P}_\kappa$ the set of all
probability measures $\pi$ on $\mathcal{X}\times\{-1,1\}$ satisfying
$\operatorname{MA}(\kappa$).

\begin{defi}\label{defoptimality}
Let $\phi$ be a loss function. The
remainder term $\gamma(n,M,\kappa,\mathcal{F},\pi)$ is called an \textbf{optimal rate
of model selection type aggregation (MS-aggregation) for the $\phi$-risk}
if the two following inequalities hold:
\begin{enumerate}[(i)]
\item[(i)] $\forall\mathcal{F}=\{f_1,\ldots,f_M\}$, there exists a statistic
$\tilde{f}_n$, depending on $\mathcal{F}$, such that $\forall\pi\in
\mathcal{P}_{\kappa}$, $\forall
n\geq1,$
%
\begin{equation}\label{DefOracle}\mathbb{E}\bigl[A^{(\phi)}(\tilde
{f}_n)-A^{(\phi)*} \bigr]\leq
\min_{f\in\mathcal{F}}\bigl(A^{(\phi)}(f)-A^{(\phi)*}
\bigr)+C_1\gamma(n,M,\kappa,\mathcal{F},\pi);
\end{equation}
\item[(ii)] $\exists\mathcal{F}=\{f_1,\ldots,f_M\}$ such that for any
statistic $\bar{f}_n$,
$\exists\pi\in\mathcal{P}_\kappa$, $\forall
n\geq1$
%
\begin{equation}\label{DefLower}\mathbb{E}\bigl[A^{(\phi)}(\bar
{f}_n)-A^{(\phi)*} \bigr]\geq
\min_{f\in\mathcal{F}}\bigl(A^{(\phi)}(f)-A^{(\phi)*}
\bigr)+C_2\gamma(n,M,\kappa,\mathcal{F},\pi).
\end{equation}
\end{enumerate}
Here, $C_1$ and $C_2$ are positive constants which may depend on
$\kappa$.
Moreover, when these \mbox{two} inequalities are satisfied,
we say that the procedure $\tilde{f}_n$, appearing in (\ref
{DefOracle}), is an \textbf{optimal\break MS-aggregate
for the $\phi$-risk}. If $\mathcal{C}$ denotes the convex hull of
$\mathcal{F}$ and if (\ref{DefOracle}) and (\ref{DefLower}) are
satisfied with
$\min_{f\in\mathcal{F}}(A^{(\phi)}(f)-A^{(\phi)*}
)$ replaced by $\min_{f\in\mathcal{C}}(A^{(\phi
)}(f)-A^{(\phi)*}
)$, then we say
that $\gamma(n,M,\kappa,\mathcal{F},\pi)$ is an \textbf{optimal rate
of convex aggregation type for the} $\phi$-\textbf{risk} and $\tilde{f}_n$ is an
\textbf{optimal convex aggregation procedure for the} $\phi$-\textbf{risk}.
\end{defi}

In Tsybakov \cite{tsy03}, the optimal rate of aggregation depends only on
$M$ and $n$. In our case, the residual term may be a function of
the underlying probability measure $\pi$, of the class $\mathcal{F}$ and
of the margin parameter $\kappa$. Note that, without any margin
assumption, we obtain $\sqrt{(\log M)/n}$ for the residual, which is
free from $\pi$ and $\mathcal{F}$. Under the margin assumption, we obtain a
residual term dependent of $\pi$ and $\mathcal{F}$ and it should be
interpreted as a normalizing factor in the ratio
\[
\frac{\mathbb{E}[A^{(\phi)}(\bar{f}_n)-A^{(\phi)*} ]-
\min_{f\in\mathcal{F}}(A^{(\phi)}(f)-A^{(\phi)*}
)}{\gamma(n,M,\kappa,\mathcal{F},\pi)}.
\]
In that case, our definition does not imply the uniqueness of
the residual.

\begin{rem} Observe that a linear function achieves its maximum over a convex
polygon at one of the vertices of the polygon. The hinge loss is
linear on $[-1,1]$ and $\mathcal{C}$ is a convex set, thus MS-aggregation
or convex aggregation of functions with values in $[-1,1]$ are
identical problems when we use the hinge loss. That is, we have
%
\begin{equation}\label{HingLinearHingeConvex}
\min_{f\in\mathcal{F}}A(f)=\min_{f\in\mathcal{C}}A(f).
\end{equation}
\end{rem}

\section{Optimal rates of convex aggregation for the hinge risk}\label{sectionoptimalhinge}
Take $M$ functions $f_1,\ldots,f_M$ with values in $[-1,1]$.
Consider the convex hull $\mathcal{C}=\operatorname{Conv}(f_1,\ldots,f_M)$. We want
to mimic the best function in $\mathcal{C}$ using the hinge risk and
working under the margin assumption. We first introduce a margin
assumption w.r.t. the hinge loss.

\begin{mah*}[Margin (or low noise) assumption for hinge risk]
The probability distribution $\pi$ on the space $\mathcal{X}\times\{-1,1
\}$ satisfies the margin assumption for hinge risk $\operatorname{MAH}(\kappa)$ with
parameter $1\leq\kappa <+\infty$ if there exists $c>0$ such
that
%
\begin{equation}\label{MAHinge}
\mathbb{E}[|f(X)-f^*(X)| ]\leq c\bigl(A(f)-A^*
\bigr)^{1/\kappa}
\end{equation}
for any function $f$ on $\mathcal{X}$ with values in $[-1,1]$.
\end{mah*}
\begin{pro}\label{proporac2}
The assumption $\operatorname{MAH}(\kappa$) is equivalent to the margin
assumption $\operatorname{MA}(\kappa$).
\end{pro}

In what follows, we will assume that $\operatorname{MA}(\kappa$) holds and thus
also that  $\operatorname{MAH}(\kappa$) holds.

The AEW aggregate of $M$ functions $f_1,\ldots,f_M$ with values in
$[-1,1]$, introduced in (\ref{AEW}) for a general loss, has a
simple form for the case of the hinge loss, given by
\begin{eqnarray}\label{aggregate}
\nonumber&&\tilde{f_n}=\sum_{j=1}^{M}w^{(n)}(f_j) f_j,\\[-8pt]\\[-8pt]
\nonumber&&\qquad\mbox{where } w^{(n)}(f_j)=\frac{\exp(
\sum_{i=1}^nY_if_j(X_i))}{\sum_{k=1}^M\exp(
\sum_{i=1}^nY_if_k(X_i))}\ \forall j=1,\ldots,M.
\end{eqnarray}

In Theorems \ref{oracleA} and \ref{optimalkappa}, we state the
optimality of our aggregates in the sense of Definition
\ref{defoptimality}.

\begin{theo}[(Oracle inequality)]\label{oracleA} Let $\kappa\geq
1$. We assume that $\pi$
satisfies $\operatorname{MA}(\kappa$). We denote by $\mathcal{C}$ the convex hull of a
finite set $\mathcal{F}$ of functions $f_1,\ldots,f_M$ with values in
$[-1,1]$. Let $\tilde{f}_n$ be either of the four aggregates
introduced in Section \ref{subsectionAggregationProcedures}. Then,
for any integers $M\geq3,n\geq1$, $\tilde{f}_n$ satisfies the
inequality
\begin{eqnarray*}
\mathbb{E}[A(\tilde{f}_n)-A^* ]&\leq&
\min_{f\in\mathcal{C}}\bigl(A(f)-A^*\bigr)\\
&&{}+C\biggl(\sqrt{\frac{\min_{f\in
\mathcal{C}}(A(f)-A^*)^{{1}/{\kappa}}\log
M}{n}}+\biggl(\frac{\log
M}{n}\biggr)^{{\kappa}/{(2\kappa-1)}}\biggr),
\end{eqnarray*}
where
$C=32(6\vee537c\vee16(2c+1/3))$ for the ERM, AERM and AEW
aggregates with $\kappa\geq1$, $c>0$ is the constant in
(\ref{MAHinge}) and $C=32(6\vee537c\vee16(2c+1/3))(2\vee
(2\kappa-1)/(\kappa-1)$ for the CAEW aggregate with $\kappa>1$.
For $\kappa=1$, the CAEW aggregate satisfies
\begin{eqnarray*}
\mathbb{E}\bigl[A\bigl(\tilde{f}^{(\mathrm{CAEW})}_n\bigr)-A^* \bigr]&\leq&
\min_{f\in\mathcal{C}}\bigl(A(f)-A^*\bigr)\\
&&{}+2C\biggl(\sqrt{\frac{\min_{f\in
\mathcal{C}}(A(f)-A^*)\log
M}{n}}+\frac{(\log M)\log n}{n}\biggr).
\end{eqnarray*}
\end{theo}

\begin{theo}[(Lower bound)]\label{optimalkappa} Let
$\kappa\geq1$ and let $M,n$ be two integers such that $2\log_2 M\leq n$. We
assume that the input space $\mathcal{X}$ is infinite. There exists an
absolute constant $C>0$, depending only on $\kappa$ and $c$, and a
set of prediction rules $\mathcal{F}=\{f_1,\ldots,f_M\}$ such that
for any
real-valued procedure $\bar{f}_n$, there exists a probability
measure $\pi$ satisfying $\operatorname{MA}(\kappa$), for which
\begin{eqnarray*}
\mathbb{E}[A(\bar{f}_n)-A^* ]&\geq& \min_{f\in\mathcal
{C}}\bigl(A(f)-A^*\bigr)
\\
&&{} +C\biggl(\sqrt{\frac{(\min_{f\in\mathcal
{C}}A(f)-A^*)^{1/\kappa}\log
M}{n}}+\biggl(\frac{\log
M}{n}\biggr)^{{\kappa}/{(2\kappa-1)}}\biggr),
\end{eqnarray*}
where
$C=c^\kappa(4\mathrm{e})^{-1}2^{-2\kappa(\kappa-1)/(2\kappa-1)}(\log
2)^{-\kappa/(2\kappa-1)}$ and $c>0$ is the constant in
(\ref{MAHinge}).
\end{theo}

Combining the exact oracle inequality of Theorem \ref{oracleA} and
the lower bound of Theorem \ref{optimalkappa}, we see that the
residual
%
\begin{equation}\label{RateOfAggregation}\sqrt{\frac{(\min_{f\in
\mathcal{C}}A(f)-A^*)^{{1}/{\kappa}}\log
M}{n}}+\biggl(\frac{\log
M}{n}\biggr)^{{\kappa}/{(2\kappa-1)}}
\end{equation}
is an
optimal rate of convex aggregation of $M$ functions with values in
$[-1,1]$ for the hinge loss. Moreover, for any real-valued
function $f$, we have $\max(1-y\psi(f(x)),0)\leq\max(1-yf(x),0)$
for all $y\in\{-1,1\}$ and $x\in\mathcal{X}$, thus
%
\begin{equation}\label{EquaProjHingeRisk}
A(\psi(f))-A^*\leq A(f)-A^*,\qquad \mbox{where }
\psi(x)=\max\bigl(-1,\min(x,1)\bigr),\ \forall x\in\mathbb{R}.
\end{equation}
Thus, by aggregating $\psi(f_1),\ldots,\psi(f_M)$, it
is easy to check that
\[
\sqrt{\frac{(\min_{f\in\mathcal{F}}A(\psi(f))-A^*)^{
{1}/{\kappa}}\log
M}{n}}+\biggl(\frac{\log M}{n}\biggr)^{{\kappa}/{(2\kappa-1)}},
\]
is an optimal rate of model-selection aggregation of $M$
 real-valued functions $f_1,\ldots,f_M$ w.r.t. the hinge loss. In both
cases, the aggregate with exponential weights, as well as ERM and
AERM, attains these optimal rates and the CAEW aggregate attains the
optimal rate if $\kappa>1$. Applications and learning properties
of the AEW procedure can be found in Lecu{\'e} \cite{lec406,lec05}
(in particular, adaptive SVM classifiers are
constructed by aggregating only $(\log n)^2$ SVM estimators). In
Theorem \ref{oracleA}, the AEW procedure satisfies an exact oracle
inequality with an optimal residual term whereas in Lecu{\'e} \cite{lec05}
and Lecu{\'e} \cite{lec406} the oracle inequalities satisfied by the AEW
procedure are not exact (there is a multiplying factor greater
than $1$ in front of the bias term) and in Lecu{\'e} \cite{lec05}, the
residual is not optimal. In Lecu{\'e} \cite{lec406}, it is proved that for
any finite set $\mathcal{F}$ of functions $f_1,\ldots,f_M$ with
values in
$[-1,1]$ and any $\epsilon>0$, there exists an absolute constant
$C(\epsilon)>0$ such that, for $\mathcal{C}$ the convex hull of
$\mathcal{F}$,
%
\begin{equation}\label{EquaLecCOLT}
\mathbb{E}\bigl[A\bigl(\tilde{f}_n^{(\mathrm{AEW})}\bigr)-A^* \bigr]\leq
(1+\epsilon) \min_{f\in\mathcal{C}}\bigl(A(f)-A^*\bigr)+C(\epsilon)
\biggl(\frac{\log
M}{n}\biggr)^{{\kappa}/{(2\kappa-1)}}.
\end{equation}
This oracle
inequality is good enough for several applications (see the
examples in Lecu{\'e} \cite{lec406}). Nevertheless, (\ref{EquaLecCOLT}) can
be easily deduced from Theorem \ref{oracleA} using Lemma
\ref{lemtsysara} and may be inefficient for constructing adaptive
estimators with exact constants (because of the factor greater than
$1$ in front of $\min_{f\in\mathcal{C}}(A(f)-A^*)$). Moreover, oracle
inequalities with a factor greater than $1$ in front of the oracle
$\min_{f\in\mathcal{C}}(A(f)-A^*)$ do not characterize the real behavior
of the  technique of aggregation which we are using. For instance, for any
strictly convex loss $\phi$, the ERM procedure satisfies (cf. Chesneau
and Lecu\'{e}
\cite{cl06})
%
\begin{equation}\label{EquaNonExcatOraclPhi}
\mathbb{E}\bigl[A^{(\phi)}\bigl(\tilde{f}_n^{(\mathrm{ERM})}\bigr)-A^{(\phi)*}
\bigr]\leq(1+\epsilon)
\min_{f\in\mathcal{F}}\bigl(A^{(\phi)}(f)-A^{(\phi)*}\bigr)+C(\epsilon
)\frac{\log
M}{n}.
\end{equation}
But, it has been recently proven, in Lecu\'{e}
\cite{lec606}, that the ERM procedure cannot mimic the oracle
faster than $\sqrt{(\log M)/n}$, whereas, for strictly convex
losses, the CAEW procedure can mimic the oracle at the rate $(\log
M)/n$ (cf. Juditsky \textit{et al.} \cite{jrt06}). Thus, for strictly convex losses, it is
better to use the aggregation procedure with exponential weights than
ERM (or even penalized ERM procedures (cf. Lecu\'{e} \cite{lec606})) to
mimic the oracle. Non-exact oracle inequalities of the form
(\ref{EquaNonExcatOraclPhi}) cannot tell us which procedure is
better to use since both ERM and CAEW procedures satisfy this
inequality.

It is interesting to note that the rate of aggregation
(\ref{RateOfAggregation}) depends on both the class $\mathcal{F}$ and
$\pi$
through the term $\min_{f\in\mathcal{C}}A(f)-A^*$. This is different from
the regression problem (cf. Tsybakov \cite{tsy03}), where the optimal
aggregation rates depend only on $M$ and $n$. Three cases can be
considered, where $\mathcal{M}(\mathcal{F},\pi)$ denotes $\min
_{f\in\mathcal{C}}(A(f)-A^*)$
and $M$ may depend on $n$ (i.e., for function classes $\mathcal{F}$
depending on $n$):
\begin{enumerate}
\item If $\mathcal{M}(\mathcal{F},\pi)\leq a (\frac{\log
M}{n})^{{\kappa}/{(2\kappa-1)}}$, for an absolute constant
$a>0$, then the hinge risk of our aggregates
attains $\min_{f\in\mathcal{C}}A(f)-A^*$ with the rate $(\frac
{\log
M}{n})^{{\kappa}/{(2\kappa-1)}}$, which can be $\log M/n$
in the case $k=1$;
\item If $ a(\frac{\log
M}{n})^{{\kappa}/{(2\kappa-1)}}\leq\mathcal{M}(\mathcal
{F},\pi)\leq
b$ for some constants $a,b>0$, then our aggregates mimic the best
prediction rule
in
$\mathcal{C}$ with a rate slower than $(\frac{\log
M}{n})^{{\kappa}/{(2\kappa-1)}}$, but faster than $((\log
M)/n)^{1/2}$;
\item If $\mathcal{M}(\mathcal{F},\pi)\geq a>0$, where $a>0$ is a
constant, then
the rate of aggregation is $\sqrt{\frac{\log M}{n}},$ as in the
case of no margin assumption.
\end{enumerate}
We can explain this behavior by the fact that not only $\kappa$, but
also $\min_{f\in\mathcal{C}}A(f)-A^*$, measures the difficulty of
classification. For instance, in the extreme case where
$\min_{f\in\mathcal{C}}A(f)-A^*=0$, which means that $\mathcal{C}$
contains the
Bayes rule, we have the fastest rate $(\frac{\log
M}{n})^{{\kappa}/{(2\kappa-1)}}$. In the worst cases, which
are realized when $\kappa$ tends to $\infty$ or
$\min_{f\in\mathcal{C}}(A(f)-A^*)\geq a>0$, where $a>0$ is an absolute
constant, the optimal rate of aggregation is the slow rate
$\sqrt{\frac{\log M}{n}}.$

\section{Optimal rates of MS-aggregation for the excess risk}\label
{sectionExcessRisk}
 We now provide oracle inequalities and lower bounds for the excess
Bayes risk. First, we can deduce, from Theorem \ref{oracleA} and
\ref{optimalkappa}, `almost optimal rates of aggregation' for the
excess Bayes risk achieved by the AEW aggregate. Second, using the
ERM aggregate, we obtain optimal rates of model selection
aggregation for the excess Bayes risk.

Using inequality (\ref{e1}), we can derive, from Theorem
\ref{oracleA}, an oracle inequality for the excess Bayes risk. The
lower bound is obtained using the same proof as in Theorem
\ref{optimalkappa}.

\begin{coro}\label{oracleR} Let $\mathcal{F}=\{f_1,\ldots,f_M\}$
be a finite set of prediction rules for an integer
$M\geq3$ and $\kappa\geq1$. We assume that $\pi$ satisfies
$\operatorname{MA}(\kappa$). Denote by $\tilde{f}_n$ either the ERM,  the AERM or
the AEW aggregate.  For any number
$a>0$ and any integer $n$, $\tilde{f}_n$ then satisfies
\begin{eqnarray}\label{approxoracleR}\mathbb{E}[R(\tilde
{f}_n)-R^* ]&\leq&
2(1+a)\min_{j=1,\ldots,M}\bigl(R(f_j)-R^*\bigr)\nonumber\\[-8pt]\\[-8pt]
\nonumber &&{}+\big[C+(C^{2\kappa
}/a)^{1/(2\kappa-1)}\big]
\biggl(\frac{\log
M}{n}\biggr)^{{\kappa}/{(2\kappa-1)}},
\end{eqnarray}
where
$C=32(6\vee537c\vee16(2c+1/3))$. The CAEW aggregate satisfies the
same inequality with $C=32(6\vee537c\vee16(2c+1/3))(2\vee
(2\kappa-1)/(\kappa-1)$ when $\kappa>1$. For $\kappa=1$, the CAEW
aggregate satisfies (\ref{approxoracleR}), where we need to multiply
the residual by $\log n$.

Moreover, there exists a finite set of prediction rules
$\mathcal{F}=\{f_1,\ldots,f_M\}$ such that, for any classifier $\bar{f}_n$,
there exists a probability measure $\pi$ on $\mathcal{X}\times\{
-1,1\}$
satisfying $\operatorname{MA}(\kappa$), such that, for any $n\geq1,a>0$,
\[
\mathbb{E}[R(\bar{f}_n)-R^* ]\geq
2(1+a)\min_{f\in\mathcal{F}}\bigl(R(f)-R^*
\bigr)+C(a)\biggl(\frac{\log M}{n}
\biggr)^{{\kappa}/{(2\kappa-1)}},
\]
where $C(a)>0$
is a constant depending only on $a$.
\end{coro}

Due to Corollary \ref{oracleR},
\[
\biggl(\frac{\log M}{n}\biggr)^{{\kappa}/{(2\kappa-1)}}
\]
is an almost
optimal rate of MS-aggregation for the excess risk and the AEW
aggregate achieves this rate. The word ``almost'' is used here because
$\min_{f\in\mathcal{F}}(R(f)-R^*)$ is multiplied by a constant
greater than $1$. Oracle inequality (\ref{approxoracleR}) is not
exact since the minimal excess risk over $\mathcal{F}$ is multiplied
by the
constant $2(1+a)>1$. This is not the case when using the ERM
aggregate, as explained in the following theorem.
\begin{theo}\label{theooracleR} Let
$\kappa\geq1$. We assume that $\pi$ satisfies $\operatorname{MA}(\kappa$). We
denote by $\mathcal{F}=\{f_1,\ldots,f_M\}$ a set of prediction rules. The
ERM aggregate over $\mathcal{F}$ satisfies, for any integer $n\geq1$,
\begin{eqnarray*}
\mathbb{E}\bigl[R\bigl(\tilde{f}^{(\mathrm{ERM})}_n\bigr)-R^* \bigr]&\leq&
\min_{f\in\mathcal{F}}\bigl(R(f)-R^*\bigr)\\
&&{}+C\biggl(\sqrt{\frac{\min_{f\in
\mathcal{F}}(R(f)-R^*)^{{1}/{\kappa}}\log
M}{n}}+\biggl(\frac{\log
M}{n}\biggr)^{\kappa/(2\kappa-1)}\biggr),
\end{eqnarray*}
where
$C=32(6\vee537c_0\vee16(2c_0+1/3))$ and $c_0$ is the constant
appearing in $\operatorname{MA}(\kappa$).
\end{theo}

Using Lemma \ref{lemtsysara}, we can deduce the results of
Herbei and Wegkamp \cite{hw05} from Theorem \ref{theooracleR}. Oracle inequalities
under $\operatorname{MA}(\kappa$) have already been stated in Massart \cite{m04} (cf.
Boucheron \textit{et al.} \cite{bbl05}), but the remainder term obtained is worse than the
one obtained in Theorem \ref{theooracleR}.

According to Definition \ref{defoptimality}, combining Theorem
\ref{theooracleR} and the following theorem, the rate
\[
\sqrt{\frac{\min_{f\in\mathcal{F}}(R(f)-R^*)^{{1}/{\kappa
}}\log
M}{n}}+\biggl(\frac{\log M}{n}\biggr)^{\kappa/(2\kappa-1)}
\]
is an optimal rate of MS-aggregation w.r.t. the excess Bayes risk.
The ERM aggregate achieves this rate.

\begin{theo}[(Lower bound)]\label{optimalkappaR} Let $M\geq3$ and
$n$ be two integers such that $2\log_2 M\leq n$ and $\kappa\geq1$.
Assume that $\mathcal{X}$ is infinite. There exists an absolute constant
$C>0$ and a set of prediction rules $\mathcal{F}=\{f_1,\ldots,f_M\}$ such
that for any procedure $\bar{f}_n$ with values in $\mathbb{R}$,
there exists a probability measure $\pi$ satisfying $\operatorname{MA}(\kappa$), for
which
\begin{eqnarray*}
\mathbb{E}[R(\bar{f}_n)-R^* ]&\geq&\min_{f\in\mathcal
{F}}\bigl(R(f)-R^*\bigr)\\
&&+{}C\biggl(\sqrt{\frac{(\min_{f\in\mathcal{F}}R(f)-R^*)^{
{1}/{\kappa}}\log
M}{n}}+\biggl(\frac{\log
M}{n}\biggr)^{\kappa/(2\kappa-1)}\biggr),
\end{eqnarray*}
where
$C={c_0}^\kappa(4\mathrm{e})^{-1}2^{-2\kappa(\kappa-1)/(2\kappa-1)}(\log
2)^{-\kappa/(2\kappa-1)}$ and $c_0$ is the constant appearing in
$\operatorname{MA}(\kappa$).
\end{theo}

\section{Proofs}\label{proofs}
\begin{pf*}{Proof of Proposition \protect\ref{proporac2}} Since, for any function
$f$ from $\mathcal{X}$ to $\{-1,1\}$, we have $2(R(f)-R^*)=A(f)-A^*$, it
follows that
$\operatorname{MA}(\kappa$) is implied by MAH($\kappa$).

Assume that $\operatorname{MA}(\kappa$) holds. We first explore the case
$\kappa>1$, where $\operatorname{MA}(\kappa$) implies that there exists a constant
$c_1>0$ such that $\mathbb{P}(|2\eta(X)-1|\leq t )\leq
c_1 t^{1/(\kappa-1)}$ for any $t>0$ (cf. Boucheron \textit{et al.} \cite{bbl05}). Let $f$ be
a function from $\mathcal{X}$ to $[-1,1]$. We have, for any $t>0$,
\begin{eqnarray*}
A(f)-A^* &=& \mathbb{E}[|2\eta(X)-1||f(X)-f^*(X)| ]
\\
&\geq&
t\mathbb{E}\bigl[|f(X)-f^*(X)|\1_{|2\eta(X)-1|\geq t}
\bigr]\\
&\geq& t\bigl(\mathbb{E}[|f(X)-f^*(X)|]-2\mathbb{P}
\bigl(|2\eta(X)-1|\leq t \bigr)
\bigr)\\
&\geq&
t\bigl(\mathbb{E}[|f(X)-f^*(X)|]-2c_1t^{1/(\kappa-1)}
\bigr).
\end{eqnarray*}
For
$t_0=((\kappa-1)/(2c_1\kappa))^{\kappa-1}\mathbb{E}[|f(X)-f^*(X)|
]^{\kappa-1}$, we obtain
\[
A(f)-A^*\geq\bigl((\kappa-1)/(2c_1\kappa)\bigr)^{\kappa-1}\kappa^{-1}
\mathbb{E}[|f(X)-f^*(X)|]^{\kappa}.
\]

For the case $\kappa=1$, $\operatorname{MA}(1$) implies that there exists $h>0$
such that $|2\eta(X)-1|\geq h$ a.s. Indeed, if for any
$N\in\mathbb{N}^*$ (the set of all positive integers), there exists
$A_N\in\mathcal{A}$ (the $\sigma$-algebra on $\mathcal{X}$) such
that $P^X(A_N)>0$
and $|2\eta(x)-1|\leq N^{-1}, \forall x\in A_N$, then, for
\[
f_N(x)=\cases{
 -f^*(x), &\quad  if
$x\in A_N$,\cr
f^*(x), &\quad  otherwise,
}
\]
we obtain
$R(f_N)-R^*\leq2P^X(A_N)/N$ and
$\mathbb{E}[|f_N(X)-f^*(X)|]=2P^X(A_N),$ and there is
no constant $c_0>0$ such that $P^X(A_N)\leq c_0 P^X(A_N)/N$ for
all $N\in\mathbb{N}^*$. So, assumption $\operatorname{MA}(1$) does not hold if no
$h>0$ satisfies $|2\eta(X)-1|\geq h$ a.s. Thus, for any $f$ from
$\mathcal{X}$ to $[-1,1]$, we have $A(f)-A^* =
\mathbb{E}[|2\eta(X)-1||f(X)-f^*(X)| ]\geq h
\mathbb{E}[|f(X)-f^*(X)| ].$
\end{pf*}
\begin{pf*}{Proof of Theorem \protect\ref{oracleA}} We   start with a general
result which says that if $\phi$ is a convex loss, then the
aggregation procedures with the weights $w^{(n)}(f),f\in\mathcal{F}$,
introduced in (\ref{AEW}) satisfy
%
\begin{equation}\label{AgregatVsERM1}A_n^{(\phi)}\bigl(\tilde
{f}_n^{(\mathrm{AEW})}\bigr)\leq
A_n^{(\phi)}\bigl(\tilde{f}_n^{(\mathrm{ERM})}\bigr)+\frac{\log M}{n}
\quad\mbox{and}\quad
A_n^{(\phi)}\bigl(\tilde{f}_n^{(\mathrm{AERM})}\bigr)\leq
A_n^{(\phi)}\bigl(\tilde{f}_n^{(\mathrm{ERM})}\bigr).
\end{equation}
Indeed, take
$\phi$ to be a convex loss. We have
$\phi(Y\tilde{f}_n(X))\leq\sum_{f\in\mathcal{F}}w^{(n)}(f)\phi(Yf(X)),$
thus
\[
A_n^{(\phi)}(\tilde{f_n})\leq\sum_{f\in\mathcal{F}}w^{(n)}(f)
A_n^{(\phi)}(f).
\]
Any $f\in\mathcal{F}$ satisfies
\[
A_n^{(\phi)}(f)=A_n^{(\phi)}\bigl(\tilde{f}^{(\mathrm{ERM})}_n\bigr)+n^{-1} \bigl(
\log\bigl(w^{(n)}\bigl(\tilde{f}^{(\mathrm{ERM})}_n\bigr)\bigr)-\log\bigl(w^{(n)}(f)\bigr)\bigr),
\]
thus,
by averaging this equality over the $w^{(n)}(f)$ and using
$\sum_{f\in\mathcal{F}}w^{(n)}(f)\log(\frac{w^{(n)}(f)}{M^{-1}}
)=K(w|u) \geq0$, where $K(w|u)$ denotes the Kullback--Leibler
divergence between the weights $w=(w^{(n)}(f))_{f\in\mathcal{F}}$ and the
uniform weights $u=(1/M)_{f\in\mathcal{F}}$, we obtain the first inequality
of~(\ref{AgregatVsERM1}). Using the convexity of $\phi$, we obtain a
similar result for the AERM aggregate.

Let $\tilde{f}_n$ be either the ERM,  the AERM or the AEW aggregate
for the class $\mathcal{F}=\{f_1,\ldots,f_M\}$. In all cases, we have, according to
(\ref{AgregatVsERM1}),
%
\begin{equation}\label{AgregatVsERM}
A_n(\tilde{f_n})\leq\min_{i=1,\ldots,M}A_n(f_i)+\frac{\log M}{n}.
\end{equation}

Let $\epsilon>0$. We consider
$\mathcal{D}=\{f\in\mathcal{C}\dvtx A(f)>A_\mathcal{C}+2\epsilon
\}$, where
$A_\mathcal{C}\stackrel{\mathrm{def}}{=}\min_{f\in\mathcal
{C}}A(f)$. Let $x>0$. If
\[
\sup_{f\in\mathcal{D}}\frac
{A(f)-A^*-(A_n(f)-A_n(f^*))}{A(f)-A^*+x}\leq
\frac{\epsilon}{A_\mathcal{C}-A^*+2\epsilon+x}
\]
then, for any $f\in\mathcal{D}$,
we have
\[
A_n(f)-A_n(f^*)\geq
A(f)-A^*-\frac{\epsilon(A(f)-A^*+x)}{A_\mathcal{C}-A^*+2\epsilon
+x}\geq
A_\mathcal{C}-A^*+\epsilon,
\]
because $A(f)-A^*\geq A_\mathcal{C}-A^*+2\epsilon$.
Hence,
%
\begin{eqnarray}\label{EquaPremIneg}
&&\mathbb{P}
\biggl[\inf_{f\in\mathcal{D}}\bigl(A_n(f)-A_n(f^*)
\bigr)<A_\mathcal{C}-A^*+\epsilon
\biggr]\\ &&\quad \leq
\mathbb{P}\biggl[\sup_{f\in\mathcal{D}}\frac
{A(f)-A^*-(A_n(f)-A_n(f^*))}{A(f)-A^*+x}>
\frac{\epsilon}{A_\mathcal{C}-A^*+2\epsilon+x}
\biggr]\nonumber.
\end{eqnarray}

According to (\ref{HingLinearHingeConvex}), for
$f'\in\{f_1,\ldots,f_M\}$ such that
$A(f')=\min_{j=1,\ldots,M}A(f_j)$, we have
$A_{\mathcal{C}}=\inf_{f\in\mathcal{C}}A(f)=\inf_{f\in\{
f_1,..,f_M\}}A(f)=A(f')$.
According to (\ref{AgregatVsERM}), we have
\[
A_n(\tilde{f_n})\leq
\min_{j=1,\ldots,M}A_n(f_j)+\frac{\log M}{n} \leq A_n(f')+\frac
{\log
M}{n}.
\]
Thus, if we assume that $A(\tilde{f_n})>A_{\mathcal{C}}+2\epsilon$,
then, by definition, we have $\tilde{f}_n\in\mathcal{D}$ and thus
there exists
$f\in\mathcal{D}$ such that $A_n(f)-A_n(f^*)\leq
A_n(f')-A_n(f^*)+(\log
M)/n$. According to (\ref{EquaPremIneg}), we have
\begin{eqnarray*}
&&\mathbb{P}[A(\tilde{f_n})>A_{\mathcal{C}}+2\epsilon ]
\\
&&\quad  \leq
\mathbb{P}\biggl[\inf_{f\in\mathcal{D}}A_n(f)-A_n(f^*)
\leq
A_n(f')-A_n(f^*)+\frac{\log M}{n} \biggr]
\\
&&\quad \leq \mathbb{P}\biggl[\inf_{f\in\mathcal{D}}
A_n(f)-A_n(f^*)\leq
A_\mathcal{C}-A^*+\epsilon\biggr]
\\
&&\qquad{}+\mathbb{P}\biggl[A_n(f')-A_n(f^*)\geq
A_\mathcal{C}-A^*+ \epsilon-\frac{\log M}{n} \biggr]
\\
&&\quad \leq
\mathbb{P}\biggl[\sup_{f\in\mathcal{C}}\frac
{A(f)-A^*-(A_n(f)-A_n(f^*))}{A(f)-A^*+x}>
\frac{\epsilon}{A_\mathcal{C}-A^*+2\epsilon+x} \biggr]\\
&&\qquad{} +\mathbb{P}\biggl[A_n(f')-A_n(f^*)\geq A_\mathcal{C}-A^*+
\epsilon-\frac{\log M}{n} \biggr].
\end{eqnarray*}
If we assume that
\[
\sup_{f\in\mathcal{C}}\frac{A(f)-A^*-(A_n(f)-A_n(f^*))}{A(f)-A^*+x}>
\frac{\epsilon}{A_\mathcal{C}-A^*+2\epsilon+x},
\]
then there exists
$f=\sum_{j=1}^{M}w_jf_j\in\mathcal{C}$ (where $w_j\geq0$ and $\sum
w_j=1$) such that
\[
\frac{A(f)-A^*-(A_n(f)-A_n(f^*))}{A(f)-A^*+x}>
\frac{\epsilon}{A_\mathcal{C}-A^*+2\epsilon+x}.
\]
The linearity of the
hinge loss on $[-1,1]$ leads to
\begin{eqnarray*}
&& \frac{A(f)-A^*-(A_n(f)-A_n(f^*))}{A(f)-A^*+x}
\\
&&\quad =
\frac{\sum_{j=1}^M
w_j[A(f_j)-A^*-(A_n(f_j)-A_n(f^*))]}{\sum_{j=1}^Mw_j[A(f_j)-A^*+x]}
\end{eqnarray*}
and, according to Lemma \ref{petitlemme}, we have
\[
\max_{j=1,\ldots,M}\frac{A(f_j)-A^*-(A_n(f_j)-A_n(f^*))}{A(f_j)-A^*+x}>
\frac{\epsilon}{A_\mathcal{C}-A^*+2\epsilon+x}.
\]

We now use the relative concentration inequality of Lemma
\ref{LemDevRela} to obtain
\begin{eqnarray*}
&&\mathbb{P}\biggl[\max_{j=1,\ldots,M}\frac
{A(f_j)-A^*-(A_n(f_j)-A_n(f^*))}{A(f_j)-A^*+x}>
\frac{\epsilon}{A_\mathcal{C}-A^*+2\epsilon+x} \biggr]\\
&&\quad \leq M\biggl(
1+\frac{8c(A_\mathcal{C}-A^*+2\epsilon+x)^2x^{1/\kappa}}{n(\epsilon
x)^2}\biggr)
\exp\biggl(-\frac{n(\epsilon x)^2}{8c(A_\mathcal{C}-A^*+2\epsilon
+x)^2x^{1/\kappa}} \biggr)\\
& &\qquad{} + M\biggl(1+\frac{16(A_\mathcal{C}-A^*+2\epsilon+x)}{3n\epsilon x}
\biggr)\exp\biggl(-\frac{3n\epsilon x}{16(A_\mathcal
{C}-A^*+2\epsilon+x)}
\biggr).
\end{eqnarray*}
Using Proposition \ref{proporac2} and Lemma
\ref{lemVarianceMargin} to upper bound the variance term and
applying Bernstein's inequality, we get
\begin{eqnarray*}
&&\mathbb{P}\biggl[A_n(f')-A_n(f^*)\geq A_\mathcal{C}-A^*+
\epsilon-\frac{\log M}{n}\biggr]\\
&&\quad\leq \exp\biggl(-\frac{n(\epsilon-(\log M)
/n)^2}{4c(A_\mathcal{C}-A^*)^{1/\kappa}+(8/3)(\epsilon-(\log M)/n)}
\biggr)
\end{eqnarray*}
for any $\epsilon>(\log M)/n$. We take
$x=A_\mathcal{C}-A^*+2\epsilon$, then, for any $(\log M)/n < \epsilon<1$,
we have
\begin{eqnarray*}
&&\mathbb{P}\bigl(A(\tilde{f}_n)>A_\mathcal{C}+2\epsilon\bigr)\\
&&\quad\leq
\exp\biggl(-\frac{n(\epsilon-\log M
/n)^2}{4c(A_\mathcal{C}-A^*)^{1/\kappa}+(8/3)(\epsilon-(\log M)/n)}
\biggr)\\
&&\qquad{}+M\biggl(1+\frac{32c(A_\mathcal{C}-A^*+2\epsilon)^{1/\kappa
}}{n\epsilon^2}
\biggr)\exp\biggl(-
\frac{n\epsilon^2}{32c(A_\mathcal{C}-A^*+2\epsilon)^{1/\kappa
}}\biggr)\\
&&\qquad{}+M\biggl(1+\frac{32}{3n\epsilon}\biggr)\exp\biggl(-\frac
{3n\epsilon}{32}
\biggr).
\end{eqnarray*}
Thus, for $2(\log M)/n<u<1$, we have
%
\begin{equation}\label{equaInit}
\mathbb{E}[A(\tilde{f}_n)-A_\mathcal{C}]\leq2u+2
\int_{u/2}^1\bigl[T_1(\epsilon)+M\bigl(T_2(\epsilon)+T_3(\epsilon)\bigr)
\bigr]\,\mathrm{d}\epsilon,
\end{equation}
where
\begin{eqnarray*}
T_1(\epsilon)&=&\exp\biggl(-\frac{n(\epsilon-(\log M)
/n)^2}{4c((A_\mathcal{C}-A^*)/2)^{1/\kappa}+(8/3)(\epsilon-(\log M)/n)}
\biggr),\\
T_2(\epsilon)&=&\biggl(1+\frac{64c(A_\mathcal{C}-A^*+2\epsilon
)^{1/\kappa}}{2^{1/\kappa}n\epsilon^2}
\biggr)\exp\biggl(-
\frac{2^{1/\kappa}n\epsilon^2}{64c(A_\mathcal{C}-A^*+2\epsilon
)^{1/\kappa}}\biggr)
\end{eqnarray*}
and
\[
T_3(\epsilon)=\biggl(1+\frac{16}{3n\epsilon}\biggr)\exp
\biggl(-\frac{3n\epsilon}{16}
\biggr).
\]

Set $\beta_1=\min(32^{-1}, (2148c)^{-1},(64(2c+1/3))^{-1})$, where
the constant $c>0$ appears in MAH($\kappa$). Consider separately
the following cases, (C$1$) and (C$2$).
\begin{enumerate}[(C1)]
\item[(C1)] \textit{The case} $A_\mathcal{C}-A^*\geq(\log
M/(\beta_1 n))^{\kappa/(2\kappa-1)}$. Denote by $\mu(M)$ the
solution of $\mu=3M\exp(-\mu)$. We have $(\log M)/2 \leq
\mu(M)\leq\log M$. Take $u$ such that
$(n\beta_1u^2)/(A_\mathcal{C}-A^*)^{1/\kappa}=\mu(M)$. Using the
definitions of case (C$1$) and   $\mu(M)$, we get $u\leq
A_\mathcal{C}-A^*$.\vadjust{\goodbreak} Moreover, $u\geq4(\log M)/n$, thus
\begin{eqnarray*}
\int_{u/2}^{1}T_1(\epsilon)\,\mathrm{d}\epsilon &\leq&
\int_{u/2}^{({A_\mathcal{C}-A^*})/{2}}\exp\biggl(-\frac
{n(\epsilon/2)^2}{(4c+4/3)(A_\mathcal{C}-A^*)^{1/\kappa}}
\biggr)\,\mathrm{d}\epsilon\\
&&{}+\int_{{(A_\mathcal{C}-A^*)}/{2}}^{1}\exp\biggl(-\frac
{n(\epsilon/2)^2}{(8c+4/3)\epsilon^{1/\kappa}}\biggr)\,\mathrm{d}\epsilon.
\end{eqnarray*}
Using Lemma \ref{intalpha} and the inequality $u\leq A_\mathcal{C}-A^*$,
we obtain
%
\begin{eqnarray}\label{equa11}
\int_{u/2}^{1}T_1(\epsilon)\,\mathrm{d}\epsilon
&\leq&
\frac{64(2c+1/3)(A_\mathcal{C}-A^*)^{1/\kappa}}{nu}\nonumber\\[-8pt]\\[-8pt]\nonumber
&&{}\times\exp\biggl(-
\frac{nu^2}{64(2c+1/3)(A_\mathcal{C}-A^*)^{1/\kappa}}\biggr).
\end{eqnarray}

 We have $128c(A_\mathcal{C}-A^*+u)\leq nu^2$. Thus, using Lemma
\ref{intalpha}, we get
%
\begin{eqnarray}\label{equa12}
\nonumber\int_{u/2}^1T_2(\epsilon)\,\mathrm{d}\epsilon
&\leq&
2\int_{u/2}^{(A_\mathcal{C}-A^*)/2}\exp\biggl(-\frac{n\epsilon
^2}{64c(A_\mathcal{C}-A^*)^{1/\kappa}}\biggr)\,\mathrm{d}\epsilon\\
&&{}+2\int_{(A_\mathcal{C}-A^*)/2}^{1}\exp\biggl(-\frac{n\epsilon
^{2-1/\kappa}}{128c}
\biggr)\,\mathrm{d}\epsilon\\
\nonumber&\leq&\frac{2148c(A_\mathcal{C}-A^*)^{1/\kappa}}{nu}\exp
\biggl(-\frac{nu^2}{2148c(A_\mathcal{C}-A^*)^{1/\kappa}}\biggr).
\end{eqnarray}
We have $u\geq32(3n)^{-1}$, so
%
\begin{eqnarray}\label{equa13}\int_{u/2}^{1}T_3(\epsilon)\,\mathrm{d}\epsilon &\leq& \frac
{64}{3n}\exp\biggl(-\frac{3nu}{64} \biggr)
\nonumber
\\[-8pt]
\\[-8pt]
\nonumber
&\leq&
\frac{64(A_\mathcal{C}-A^*)^{1/\kappa}}{3nu}\exp\biggl(-\frac
{3nu^2}{64(A_\mathcal{C}-A^*)^{1/\kappa}}\biggr).
\end{eqnarray}

 From (\ref{equa11}), (\ref{equa12}), (\ref{equa13}) and
(\ref{equaInit}), we obtain
\[
\mathbb{E}[A(\tilde{f}_n)-A_\mathcal{C}]\leq
2u+6M\frac{(A_\mathcal{C}-A^*)^{1/\kappa}}{n\beta_1u}
\exp\biggl(-\frac{n\beta_1u}{(A_\mathcal{C}-A^*)^{1/\kappa}}
\biggr).
\]
The
definitions of $u$ leads to $\mathbb{E}[A(\tilde
{f}_n)-A_\mathcal{C}]\leq
4\sqrt{\frac{(A_\mathcal{C}-A^*)^{1/\kappa}\log M}{n\beta_1}}.$

\item[(C2)] \textit{The case} $A_\mathcal{C}-A^*\leq(\log
M/(\beta_1
n))^{\kappa/(2\kappa-1)}$. We now choose $u$ such that $n\beta_2
u^{(2\kappa-1)/\kappa}=\mu(M)$, where $\beta_2=\min
(3(32(6c+1))^{-1}, (256c)^{-1}, 3/64)$. Using the definition of case
(C$2$) and  $\mu(M)$, we get $u\geq A_\mathcal{C}-A^*$. Using Lemma
\ref{intalpha} and $u>4(\log M)/n$, $u\geq
2(32c/n)^{\kappa/(2\kappa-1)}$ and $u>32/(3n)$, respectively, we
obtain
%
\begin{eqnarray}\label{equa21}
\int_{u/2}^{1}T_1(\epsilon)\,\mathrm{d}\epsilon
&\leq&
\frac{32(6c+1)}{3nu^{1-1/\kappa}}\exp\biggl(-\frac{3nu^{2-1/\kappa
}}{32(6c+1)}
\biggr),
\nonumber
\\[-8pt]
\\[-8pt]
\nonumber
\int_{u/2}^{1}T_2(\epsilon)\,\mathrm{d}\epsilon
&\leq& \frac{128c}{nu^{1-1/\kappa}}\exp\biggl(-\frac{nu^{2-1/\kappa}}{128c}
\biggr)
\end{eqnarray}
and
%
\begin{equation}\label{equa23}\int_{u/2}^{1}T_3(\epsilon)\,\mathrm{d}\epsilon
\leq
\frac{64}{3nu^{1-1/\kappa}}\exp\biggl(-\frac{3nu^{2-1/\kappa}}{64}
\biggr).
\end{equation}
From (\ref{equa21}), (\ref{equa23}) and
(\ref{equaInit}), we obtain
\[
\mathbb{E}[A(\tilde{f}_n)-A_\mathcal{C}]\leq
2u+6M\frac{\exp(-n\beta_2 u^{(2\kappa-1)/\kappa}
)}{n\beta_2 u^{1-1/\kappa}}.
\]
The definition of $u$ yields
$\mathbb{E}[A(\tilde{f}_n)-A_\mathcal{C}]\leq4
(\frac{\log
M}{n\beta_2} )^{\kappa/(2\kappa-1)}.$
\end{enumerate}
 Finally, we obtain
\[
\mathbb{E}[A(\tilde{f}_n)-A_\mathcal{C}]\leq4 \cases{\displaystyle{
\biggl(\frac{\log M}{n\beta_2} \biggr)^{\kappa/(2\kappa-1)}}, &\quad
if  $A_\mathcal{C}-A^*\leq\biggl(\dfrac{\log M}{n\beta_1}
\biggr)^{\kappa/(2\kappa-1)}$, \cr
\displaystyle{\sqrt{\frac{(A_\mathcal{C}-A^*)^{1/\kappa}\log
M}{n\beta_1}}},
&\quad
otherwise.
}
\]
%
For the CAEW aggregate, it suffices to upper bound the sums by
integrals in the following inequality to get the result:
\begin{eqnarray*}
\mathbb{E}\bigl[A\bigl(\tilde{f}_n^{(\mathrm{CAEW})}\bigr)-A^* \bigr] &\leq&
\frac{1}{n}\sum_{k=1}^n \mathbb{E}\bigl[A\bigl(\tilde{f}_k^{(\mathrm{AEW})}\bigr)-A^*
\bigr]\\
 & \leq& \min_{f\in\mathcal{C}}A(f)-A^*+
C\Biggl\{\sqrt{(A_\mathcal{C}-A^*)^{1/\kappa}\log M}
\Biggl(\frac{1}{n}\sum_{k=1}^n\frac{1}{\sqrt{k}}\Biggr)\\
&&{}\hspace*{92pt}+
(\log M)^{\kappa/(2\kappa-1)}\frac{1}{n}\sum_{k=1}^n\frac{1}
{k^{\kappa/(2\kappa-1)}} \Biggr\}.
\end{eqnarray*}\upqed
\end{pf*}

\textbf{Proof of Theorem \ref{optimalkappa}.} Let $a$ be a positive
number, $\mathcal{F}$ be a finite set of $M$ real-valued functions and
$f_1,\ldots,f_M$ be $M$ prediction rules (which will be carefully chosen
 in what follows). Using~(\ref{HingLinearHingeConvex}), taking
$\mathcal{F}=\{f_1,\ldots,f_M\}$ and assuming that
$f^*\in\{f_1,\ldots,f_M\}$, we\vadjust{\goodbreak} obtain\looseness=1
%
\begin{eqnarray}\label{IneqBornInf}
&&\inf_{\hat{f}_n}\sup_{\pi\in\mathcal{P}_\kappa}
\biggl(\mathbb{E}[A(\hat{f}_n)-A^* ]-(1+a)\min_{f\in
{\operatorname{Conv}}(\mathcal{F})}\bigl(A(f)-A^*\bigr) \biggr)
\nonumber
\\[-8pt]
\\[-8pt]
\nonumber
&&\qquad\geq
\inf_{\hat{f}_n}\mathop{\sup_{\pi\in\mathcal{P}}}_{\kappa
f^*\in\{f_1,\ldots,f_M\}} \mathbb{E}[A(\hat{f}_n)-A^*
],
\end{eqnarray}
where ${\operatorname{Conv}}(\mathcal{F})$ is the set made of
all  convex combinations of elements in $\mathcal{F}$. Let $N$ be an
integer such that $2^{N-1}\leq M$, $x_1,\ldots,x_N$ be $N$ distinct
points of $\mathcal{X}$ and $w$ be a positive number satisfying
$(N-1)w\leq1$. Denote by $P^X$ the probability measure on $\mathcal
{X}$ such
that $P^X(\{x_j\})=w$, for $j=1,\ldots,N-1$, and
$P^X(\{x_N\})=1-(N-1)w$. We consider the cube
$\Omega=\{-1,1\}^{N-1}$. Let $0<h<1$. For all
$\sigma=(\sigma_1,\ldots,\sigma_{N-1})\in\Omega$ we consider
\[
\eta_{\sigma}(x)=\cases{
 (1+\sigma_j h)/2, &\quad if  $x=x_1,\ldots
,x_{N-1}$,\cr
1, &\quad if  $x =x_N$.
}
\]
For all $\sigma\in\Omega$, we denote by $\pi_\sigma$ the probability
measure on $\mathcal{X}\times\{-1,1\}$ having $P^X$ for marginal on
$\mathcal{X}$
and $\eta_{\sigma}$ for conditional probability function.

Assume that $\kappa>1$. We have
$\mathbb{P}(|2\eta_\sigma(X)-1|\leq t )=(N-1)w\1_{h\leq
t}$ for any $0\leq t <1$. Thus, if we assume that $(N-1)w\leq
h^{1/(\kappa-1)}$, then $\mathbb{P}(|2\eta_\sigma(X)-1|\leq t
)\leq t^{1/(\kappa-1)}$ for all $0\leq t <1$. Thus, according
to Tsybakov \cite{tsy04}, $\pi_\sigma$ belongs to $\mathcal{P}_\kappa$.

We denote by $\rho$ the Hamming distance on $\Omega$. Let $\sigma,
\sigma'\in\Omega$ be such that $\rho(\sigma,\sigma')=1$. Denote by $H$
the Hellinger distance. Since $H^2(\pi_\sigma^{\otimes n},
\pi_{\sigma'}^{\otimes
n})=2(1-(1-H^2(\pi_{\sigma},\pi_{\sigma'})/2
)^n)$
and
\begin{eqnarray*}
H^2(\pi_{\sigma},\pi_{\sigma'})&=&w\sum_{j=1}^{N-1}
\bigl(\sqrt{\eta_{\sigma}(x_j)}-\sqrt{\eta_{\sigma'}(x_j)}\bigr)^2+
\bigl(\sqrt{1-\eta_{\sigma}(x_j)}-\sqrt{1-\eta_{\sigma
'}(x_j)}\bigr)^2\\
&=&2w\bigl(1-\sqrt{1-h^2}\bigr),
\end{eqnarray*}
the Hellinger distance between the measures
$\pi_\sigma^{\otimes n}$ and $\pi_{\sigma'}^{\otimes n}$ satisfies
\[
H^2(\pi_\sigma^{\otimes n}, \pi_{\sigma'}^{\otimes
n})=2\bigl(1-\bigl(1-w\bigl(1-\sqrt{1-h^2}\bigr)\bigr)^n\bigr).
\]

Take $w$ and $h$ such that $w(1-\sqrt{1-h^2})\leq n^{-1}.$ Then,
$H^2(\pi_\sigma^{\otimes n}, \pi_{\sigma'}^{\otimes
n})\leq\beta=2(1-\mathrm{e}^{-1})<2$ for any integer $n$.

Let $\sigma\in\Omega$ and $\hat{f}_n$ be an estimator with values in
$[-1,1]$ (according to (\ref{EquaProjHingeRisk}), we consider only
estimators in $[-1,1]$). Using $\operatorname{MA}(\kappa$), we have, conditionally
on the observations $D_n$ and for $\pi=\pi_\sigma$,
\[
A(\hat{f}_n)-A^* \geq \bigl(
c\mathbb{E}_{\pi_\sigma}[|\hat{f}_n(X)-f^*(X)|]
\bigr)^\kappa
\geq (cw)^\kappa\Biggl(
\sum_{j=1}^{N-1}|\hat{f}_n(x_j)-\sigma_j|\Biggr)^{\kappa}.
\]
Taking here the expectations, we find $\mathbb{E}_{\pi_\sigma}[
A(\hat{f}_n)-A^*]\geq(cw)^\kappa\mathbb{E}_{\pi_\sigma
}[
(
\sum_{j=1}^{N-1}|\hat{f}_n(x_j)-\sigma_j|)^{\kappa}].$
Using Jensen's inequality and Lemma \ref{lem2}, we obtain

\begin{equation}\label{IneqAssouadApplique}\inf_{\hat{f}_n}\sup
_{\sigma\in\Omega}
\bigl(\mathbb{E}_{\pi_\sigma}[A(\hat{f}_n)-A^*
]\bigr)\geq(cw)^\kappa\biggl(\frac{N-1}{4\mathrm{e}^2}
\biggr)^\kappa.
\end{equation}

Now take $w=(nh^2)^{-1}$, $N=\lceil\log M/ \log2 \rceil$ and
$h=(n^{-1}\lceil\log M/ \log2 \rceil
)^{(\kappa-1)/(2\kappa-1)}$. Replace $w$ and $N$ in
(\ref{IneqAssouadApplique}) by these values. Thus, from
(\ref{IneqBornInf}), there exist $f_1,\ldots,f_M$ (the first $2^{N-1}$
 are ${\rm sign}(2\eta_{\sigma}-1)$ for $\sigma\in\Omega$
and any choice is allowed for the remaining $M-2^{N-1}$) such that, for any
procedure $\bar{f}_n$, there exists a probability measure $\pi$
satisfying $\operatorname{MA}(\kappa$), such that $
\mathbb{E}[A(\hat{f}_n)-A^*
]-(1+a)\min_{j=1,\ldots,M}(A(f_j)-A^*)\geq C_0
(\frac{\log M}{n} )^{\kappa/(2\kappa-1)},$ where
$C_0=c^\kappa(4\mathrm{e})^{-1}2^{-2\kappa(\kappa-1)/(2\kappa-1)}(\log
2)^{-\kappa/(2\kappa-1)}$.

Moreover, according to Lemma \ref{lemtsysara}, we have
\begin{eqnarray*}
&& a\min_{f\in\mathcal{C}}\bigl(A(f)-A^*\bigr)+
\frac{C_0}{2}\biggl(\frac{\log M}{n}
\biggr)^{\kappa/(2\kappa-1)}
\\
&&\quad \geq
\sqrt{2^{-1}a^{1/\kappa}C_0}\sqrt{\frac{(\min_{f\in\mathcal
{C}}A(f)-A^*)^{{1}/{\kappa}}\log
M}{n}}.
\end{eqnarray*}
Thus,
\begin{eqnarray*}
\mathbb{E}[A(\hat{f}_n)-A^* ]&\geq&
\min_{f\in\mathcal{C}}\bigl(A(f)-A^*\bigr)+\frac{C_0}{2}\biggl( \frac{\log
M}{n}\biggr)^{\kappa/(2\kappa-1)}\\
&&{}+\sqrt{2^{-1}a^{1/\kappa}C_0}
\sqrt{\frac{(A_\mathcal{C}-A^*)^{{1}/{\kappa}}\log
M}{n}}.
\end{eqnarray*}

For $\kappa=1$, we take $h=1/2$. Then, $|2\eta_\sigma(X)-1|\geq
1/2$ a.s., so $\pi_\sigma\in$$\operatorname{MA}(1$). It then suffices to take
$w=4/n$ and $N=\lceil\log M/ \log2 \rceil$ to obtain the result.

\begin{pf*}{Proof of Corollary \protect\ref{oracleR}} The result follows from
Theorems \ref{oracleA} and \ref{optimalkappa}. Using inequality (\ref{e1}),  Lemma \ref{lemtsysara} and the fact that
for any prediction rule $f$, we have $A(f)-A^*=2(R(f)-R^*)$, for any $a>0$,
with $t=a(A_\mathcal{C}-A^*)$ and $v=(C^2(\log
M)/n)^{\kappa/(2\kappa-1)}a^{-1/(2\kappa-1)}$, we obtain the
result.
\end{pf*}

\begin{pf*}{Proof of Theorem \protect\ref{theooracleR}} Denote by $\tilde{f}_n$ the ERM
aggregate over $\mathcal{F}$. Let $\epsilon>0$. Denote by $\mathcal
{F}_\epsilon$ the
set $\{f\in\mathcal{F}:R(f)>R_\mathcal{F}+2\epsilon\}$, where
$R_\mathcal{F}=\min_{f\in\mathcal{F}}R(f).$

Let $x>0$. If
\[
\sup_{f\in\mathcal{F}_\epsilon}\frac
{R(f)-R^*-(R_n(f)-R_n(f^*))}{R(f)-R^*+x}\leq
\frac{\epsilon}{R_\mathcal{F}-R^*+2\epsilon},
\]
then the same argument as in
Theorem \ref{oracleA} yields $R_n(f)-R_n(f^*)\geq
R_\mathcal{F}-R^*+\epsilon$ for any $f\in\mathcal{F}_\epsilon$.
So, we have
\begin{eqnarray*}
&&\mathbb{P}\biggl[
\inf_{f\in\mathcal{F}_\epsilon}R_n(f)-R_n(f^*)<R_\mathcal
{F}-R^*+\epsilon\biggr]\\
&&\quad \leq
\mathbb{P}\biggl[\sup_{f\in\mathcal{F}_\epsilon}\frac
{R(f)-R^*-(R_n(f)-R_n(f^*))}{R(f)-R^*+x}>
\frac{\epsilon}{R_\mathcal{F}-R^*+2\epsilon+x} \biggr].
\end{eqnarray*}

We consider $f'\in\mathcal{F}$ such that $\min_{f\in\mathcal
{F}}R(f)=R(f').$ If
$R(\tilde{f}_n)>R_\mathcal{F}+2\epsilon$, then $\tilde{f}_n\in
\mathcal{F}_\epsilon$, so there exists
$g\in\mathcal{F}_\epsilon$ such that $R_n(g)\leq R_n(f')$. Hence,
using the
same argument as in Theorem \ref{oracleA}, we obtain
\begin{eqnarray*}
\mathbb{P}[R(\tilde{f}_n)>R_\mathcal{F}+2\epsilon]
&\leq&
\mathbb{P}\biggl[\sup_{f\in\mathcal{F}}\frac
{R(f)-R^*-(R_n(f)-R_n(f^*))}{R(f)-R^*+x}\geq
\frac{\epsilon}{R_\mathcal{F}-R^*+2\epsilon+x} \biggr]\\
&&{}+
\mathbb{P}[R_n(f')-R_n(f^*)>R_\mathcal{F}-R^*+\epsilon].
\end{eqnarray*}

We complete the proof by using Lemma \ref{LemDevRela}, the fact that for any $f$ from $\mathcal{X}$
to $\{-1,1\}$, we have $2(R(f)-R^*)=A(f)-A^*$,
and the same arguments as those developed at the end of the proof
of Theorem \ref{oracleA}.
\end{pf*}
\begin{pf*}{Proof of Theorem \protect\ref{optimalkappaR}} Using the same argument
as the one used in the beginning of the proof of Theorem
\ref{optimalkappa}, we have, for all prediction rules
$f_1,\ldots,f_M$ and $a>0$,
\begin{eqnarray*}
&&\sup_{g_1,\ldots,g_M}\inf_{\hat{f}_n}\sup_{\pi\in\mathcal
{P}_\kappa}
\biggl(\mathbb{E}[R(\hat{f}_n)-R^*
]-(1+a)\min_{j=1,\ldots,M}\bigl(R(g_j)-R^*\bigr) \biggr)
\\
&&\quad\geq
\inf_{\hat{f}_n}\mathop{\sup_{\pi\in\mathcal{P}_\kappa}}\limits_{f^*\in\{f_1,\ldots,f_M\}} \mathbb{E}[R(\hat{f}_n)-R^*
].
\end{eqnarray*}
Consider the set of probability measures
$\{\pi_\sigma, \sigma\in\Omega\}$ introduced in the proof of
Theorem \ref{optimalkappa}. Assume that $\kappa>1$. Since for any
$\sigma\in\Omega$ and any classifier $\hat{f}_n$, we have, by
using $\operatorname{MA}(\kappa$),
\[
\mathbb{E}_{\pi_\sigma}[ R(\hat{f}_n)-R^*]
\geq(c_0w)^\kappa\mathbb{E}_{\pi_\sigma}\Biggl[ \Biggl(
\sum_{j=1}^{N-1}|\hat{f}_n(x_j)-\sigma_j|\Biggr)^{\kappa}\Biggr],
\]
using Jensen's inequality and Lemma \ref{lem2}, we obtain
\[
\inf_{\hat{f}_n}\sup_{\sigma\in\Omega}
\bigl(\mathbb{E}_{\pi_\sigma}[R(\hat{f}_n)-R^*
]\bigr)\geq(c_0w)^\kappa\biggl(\frac{N-1}{4\mathrm{e}^2}
\biggr)^\kappa.
\]

By taking $w=(nh^2)^{-1}$, $N=\lceil\log M/ \log2 \rceil$ and
$h=(n^{-1}\lceil\log M/ \log2 \rceil
)^{{(\kappa-1)}/{(2\kappa-1)}}$, there exist $f_1,\ldots,f_M$
(the first $2^{N-1}$  are ${\rm sign}(2\eta_{\sigma}-1)$ for
$\sigma\in\Omega$ and any choice is allowed for the remaining  $M-2^{N-1}$)
such that for any procedure $\bar{f}_n$, there exists a probability
measure $\pi$ satisfying $\operatorname{MA}(\kappa$), such that $
\mathbb{E}[R(\hat{f}_n)-R^*
]-(1+a)\min_{j=1,\ldots,M}(R(f_j)-R^*)\geq C_0
(\frac{\log M}{n} )^{\kappa/(2\kappa-1)},$ where
$C_0={c_0}^{\kappa}(4\mathrm{e})^{-1}2^{-2\kappa(\kappa-1)/(2\kappa-1)}(\log
2)^{-\kappa/(2\kappa-1)}$. Moreover, according to Lemma
\ref{lemtsysara}, we have
\begin{eqnarray*}
&&a\min_{f\in\mathcal{F}}\biggl(R(f)-R^*\biggr)+
\frac{C_0}{2}\biggl(\frac{\log M}{n}
\biggr)^{\kappa/(2\kappa-1)}\\
&&\quad\geq
\sqrt{a^{1/\kappa}C_0/2}\sqrt{\frac{(\min_{f\in\mathcal
{F}}R(f)-R^*)^{{1}/{\kappa}}\log
M}{n}}.
\end{eqnarray*}
The case $\kappa=1$ is treated in the same way as in the
proof of Theorem \ref{optimalkappa}.

\begin{lem}\label{intalpha}
Let $\alpha\geq1$ and $a,b>0$. An integration by parts yields
\[
\int_a^{+\infty}\exp(-bt^\alpha)\,\mathrm{d}t\leq\frac{\exp
(-ba^\alpha)}{\alpha b
a^{\alpha-1}}.
\]
\end{lem}
\begin{lem}\label{petitlemme}
Let $b_1,\ldots,b_M$ be $M$ positive numbers and $a_1,\ldots,a_M$
some numbers. We have
\[
\frac{\sum_{j=1}^Ma_j}{\sum_{j=1}^Mb_j}\leq
\max_{j=1,\ldots,M}\biggl(\frac{a_j}{b_j} \biggr).
\]
\end{lem}
\end{pf*}
\begin{pf}
\[
\sum_{j=1}^M b_j\max_{k=1,\ldots,M}\biggl(\frac{a_k}{b_k}
\biggr)\geq\sum_{j=1}^M b_j\frac{a_j}{b_j}=\sum_{j=1}^Ma_j.
\]
\end{pf}
\begin{lem}\label{lemtsysara}
Let $v,t>0$ and $\kappa\geq1$. The concavity of the logarithm
yields
\[
t+v \geq t^{{1}/{(2\kappa)}} v^{{(2\kappa-1)}/{(2\kappa)}}.
\]
\end{lem}
%
%
\begin{lem}\label{lemVarianceMargin}
Let $f$ be a function from $\mathcal{X}$ to $[-1,1]$ and $\pi$ a
probability measure on $\mathcal{X}\times\{-1,1\}$ satisfying
$\operatorname{MA}(\kappa$) for some $\kappa\geq1$. Denote by $\mathbb{V}$ the
symbol of
variance. We have
\[
\mathbb{V}\bigl(Y\bigl(f(X)-f^*(X)\bigr)\bigr)\leq c \bigl(A(f)-A^*\bigr)^{1/\kappa}
\]
and
\[
\mathbb{V}\big(\1_{Yf(X)\leq0}-\1_{Yf^*(X)\leq0}\big)\leq
c\bigl(R(f)-R^*\bigr)^{1/\kappa}.
\]
\end{lem}
\begin{lem}\label{LemDevRela} Let $\mathcal{F}=\{f_1,\ldots,f_M\}$
be a finite set of
functions from $\mathcal{X}$ to $[-1,1]$. Assume that $\pi$ satisfies
$\operatorname{MA}(\kappa$) for some $\kappa\geq1$. We have, for any positive
numbers $t,x$ and any integer $n$,
\[
 \mathbb{P}\biggl[
\max_{f\in\mathcal{F}}Z_x(f)>t\biggr] \leq
M\biggl(\biggl(1+\frac{8cx^{1/\kappa}}{n(tx)^2} \biggr)\exp\biggl(
-\frac{n(tx)^2}{8cx^{1/\kappa}}\biggr)+\biggl(1+\frac{16}{3ntx}
\biggr)\exp\biggl(-\frac{3ntx}{16} \biggr) \biggr),
\]
where the constant $c>0$ appears in $\operatorname{MAH}(\kappa$) and
$Z_x(f)=\frac{A(f)-A_n(f)-(A(f^*)-A_n(f^*))}{A(f)-A^*+x}.$
\end{lem}
\begin{pf} For any integer $j$, consider the set
$\mathcal{F}_j=\{f\in\mathcal{F}\dvtx jx\leq A(f)-A^*<(j+1)x
\}$. Using
Bernstein's inequality, Proposition \ref{proporac2} and Lemma
\ref{lemVarianceMargin} to upper bound the variance term, we obtain
\begin{eqnarray*}
&&\mathbb{P}\biggl[\max_{f\in\mathcal{F}} Z_x(f)>t \biggr]
\\
&&\quad\leq
\sum_{j=0}^{+\infty}\mathbb{P}\biggl[\max_{f\in\mathcal{F}_j} Z_x(f)>t
\biggr]
\\
&&\quad \leq \sum_{j=0}^{+\infty}
\mathbb{P}\biggl[\max_{f\in\mathcal
{F}_j}A(f)-A_n(f)-\bigl(A(f^*)-A_n(f^*)\bigr)>t(j+1)x
\biggr]
\\
&&\quad \leq
M\sum_{j=0}^{+\infty}\exp\biggl(-\frac
{n[t(j+1)x]^2}{4c((j+1)x)^{1/\kappa}+(8/3)t(j+1)x}
\biggr)
\\
&&\quad\leq
M\Biggl(\sum_{j=0}^{+\infty}\exp\biggl(-\frac
{n(tx)^2(j+1)^{2-1/\kappa}}{8cx^{1/\kappa}}
\biggr) +\exp\biggl(-(j+1)\frac{3ntx}{16} \biggr) \Biggr)
\\
&&\quad \leq
M\biggl(\exp\biggl(-\frac{nt^2x^{2-1/\kappa}}{8c}\biggr)+\exp
\biggl(-\frac{3ntx}{16}
\biggr)\biggr)
\\
& &\qquad{}
+M\int_{1}^{+\infty}\biggl(\exp\biggl(-\frac{nt^2x^{2-1/\kappa
}}{8c}u^{2-1/\kappa}
\biggr)+\exp\biggl(-\frac{3ntx}{16}u \biggr)\biggr)\,\mathrm{d}u .
\end{eqnarray*}
Lemma \ref{intalpha} leads to the result.
\begin{lem}\label{lem2}
Let $\{P_{\omega}/\omega\in\Omega\}$ be a set of probability
measures on a measurable space $(\mathcal{X},\mathcal{A})$, indexed
by the cube
$\Omega=\{0,1\}^m$ . Denote by $\mathbb{E}_\omega$ the expectation
under $P_\omega$ and by $\rho$ the Hamming distance on $\Omega$.
Assume that
\[
\forall\omega,\omega' \in\Omega/ \rho(\omega,\omega')=1,\qquad
H^2(P_{\omega},P_{\omega'})\leq\alpha<2,
\]
Then,
\[
\inf_{\hat{w}\in[0,1]^m}\max_{\omega\in\Omega}\mathbb{E}_\omega
\Biggl[\sum_{j=1}^m|\hat{w_j}-w_j |\Biggr]\geq\frac
{m}{4}\biggl(1-\frac{\alpha}{2}\biggr)^2.
\]
\end{lem}
\end{pf}
\begin{pf} Obviously, we can replace $\inf_{\hat{w}\in[0,1]^m}$
by $(1/2)\inf_{\hat{w}\in\{0,1\}^m}$ since for all $w\in\{0,1\}$ and
$\hat{w}\in[0,1]$, there exists $\tilde{w}\in\{0,1\}$ (e.g.,
the projection of $\hat{w}$ on to $\{0,1\}$) such that $|\hat{w}-w|\geq
(1/2)|\tilde{w}-w|$. We then  use Theorem $2.10$   of
Tsybakov \cite{tsy03}, page 103.
\end{pf}

\printhistory


\begin{thebibliography}{00}
\bibitem{at05}
Audibert, J.-Y. and Tsybakov, A.B. (2007).
Fast learning rates for plug-in classifiers under margin condition.
\textit{Ann. Statist.} \textbf{35}. To appear.

\bibitem{bfls98}
Bartlett, P.L., Freund, Y., Lee, W.S. and Schapire, R.E. (1998).
Boosting the margin: A new explanation for the effectiveness of
voting methods.
\textit{Ann. Statist.} \textbf{26} 1651--1686.
\MR{1673273}

\bibitem{bjm03}
Bartlett, P.L., Jordan, M.I. and McAuliffe, J.D. (2006).
Convexity, {c}lassification and {r}isk {b}ounds.
\textit{J. Amer. Statist. Assoc.} \textbf{101} 138--156.
\MR{2268032}

\bibitem{b04}
Birg{\'e}, L. (2006).
Model selection via testing: An alternative to (penalized) maximum
likelihood estimators.
\textit{Ann. Inst. H. Poincar\'{e} Probab. Statist.} \textbf{42}
273--325.
\MR{2219712}

\bibitem{bbm04}
Blanchard, G., Bousquet, O. and Massart, P. (2004).
Statistical {p}erformance of {s}upport {v}ector {m}achines.
Available at
\href{http://mahery.math.u-psud.fr/~blanchard/publi/}{http//mahery.math.u-psud.fr/\texttildelow blanchard/publi/}.

\bibitem{blv03}
Blanchard, G., Lugosi, G. and Vayatis, N. (2003).
On the rate of convergence of regularized boosting classifiers.
\textit{J. Mach. Learn. Res.} \textbf{4} 861--894.
\MR{2076000}

\bibitem{bbl05}
Boucheron, S., Bousquet, O. and Lugosi, G. (2005).
Theory of classification: A survey of some recent advances.
\textit{ESAIM Probab. Statist.} \textbf{9} 323--375.
\MR{2182250}

\bibitem{by02}
B{\"u}hlmann, P. and Yu, B. (2002).
Analyzing bagging.
\textit{Ann. Statist.} \textbf{30} 927--961.
\MR{1926165}

\bibitem{btw105}
Bunea, F., Tsybakov, A.B. and Wegkamp, M. (2005).
Aggregation for Gaussian regression.
\textit{Ann. \mbox{Statist}.} To appear. Available at
\url{http://www.stat.fsu.edu/~wegkamp}.

\bibitem{cat99}
Catoni, O. (1999).
``Universal'' aggregation rules with exact bias bounds.
Preprint n. 510, LPMA. Available at
\url{http://www.proba.jussieu.fr/mathdoc/preprints/index.html}.

\bibitem{catbook01}
Catoni, O. (2001).
\textit{Statistical Learning Theory and Stochastic Optimization.
Ecole d'{\'E}t{\'e} de Probabilit{\'e}s de Saint-Flour 2001. Lecture
Notes in Math.} \textbf{1851}. New York: Springer.
\MR{2163920}

\bibitem{cl06}
Chesneau, C. and Lecu{\'e}, G. (2006).
Adapting to unknown smoothness by aggregation of thresholded wavelet
estimators.
Submitted.

\bibitem{cv95}
Cortes, C. and Vapnik, V. (1995).
Support-vector networks.
\textit{Machine Learning} \textbf{20} 273--297.

\bibitem{dgl96}
Devroye, L., Gy\"orfi, L. and Lugosi, G. (1996).
\textit{A Probabilistic Theory of Pattern Recognition}.
New York: Springer.
\MR{1383093}

\bibitem{fs97}
Freund, Y. and Schapire, R. (1997).
A decision-theoric generalization of on-line learning and an
application to boosting.
\textit{J. Comput. Syst. Sci.} \textbf{55}
119--139.
\MR{1473055}

\bibitem{fht00}
Friedman, J., Hastie, T. and Tibshirani, R. (2000).
Additive logistic regression: A statistical view of boosting (with discussion).
\textit{Ann. Statist.} \textbf{28} 337--407.
\MR{1790002}

\bibitem{hw05}
Herbei, R. and Wegkamp, H. (2006).
Classification with reject option.
\textit{Canad. J. Statist.} \textbf{34} 709--721.

\bibitem{jn00}
Juditsky, A. and Nemirovski, A. (2000).
Functional aggregation for nonparametric estimation.
\textit{Ann. Statist.} \textbf{28} 681--712.
\MR{1792783}

\bibitem{jrt06}
Juditsky, A., Rigollet, P. and Tsybakov, A.B. (2006).
Learning by mirror averaging.
Preprint n. 1034, Laboratoire de Probabilit{\'e}s et Mod{\`e}le
al{\'e}atoires, Univ. Paris 6 and Paris 7. Available at
\url{http://www.proba.jussieu.fr/mathdoc/preprints/index.html\#2005}.

\bibitem{lec406}
Lecu{\'e}, G. (2006).
Optimal oracle inequality for aggregation of classifiers under low
noise condition.
In \textit{Proceedings of the 19th Annual Conference on Learning Theory,
COLT 2006} \textbf{32} 364--378.
\MR{2280618}

\bibitem{lec05}
Lecu{\'e}, G. (2007).
Simultaneous adaptation to the margin and to complexity in
classification.
\textit{Ann. Statist.} To appear. Available at
\url{http://hal.ccsd.cnrs.fr/ccsd-00009241/en/}.

\bibitem{lec606}
Lecu{\'e}, G. (2007).
{S}uboptimality of {p}enalized {e}mpirical {r}isk {m}inimization.
In \textit{COLT07}. To appear.

\bibitem{l99}
Lin, Y. (1999).
A note on margin-based loss functions in classification.
Technical Report 1029r,
Dept. Statistics,
Univ. Wisconsin, Madison.

\bibitem{lv04}
Lugosi, G. and Vayatis, N. (2004).
On the {B}ayes-risk consistency of regularized boosting methods.
\textit{Ann. Statist.} \textbf{32} 30--55.
\MR{2051000}

\bibitem{mt99}
Mammen, E. and Tsybakov, A.B. (1999).
Smooth discrimination analysis.
\textit{Ann. Statist.} \textbf{27} 1808--1829.
\MR{1765618}

\bibitem{m00}
Massart, P. (2000).
Some applications of concentration inequalities to statistics.
\textit{Ann. Fac. Sci. Toulouse Math.} (\textit{6}) \textbf{2} 245--303.
\MR{1813803}

\bibitem{m04}
Massart, P. (2004).
Concentration inequalities and model selection.
\textit{Lectures Notes of Saint Flour}.

\bibitem{mn03}
Massart, P. and N{\'e}d{\'e}lec, E. (2006).
Risk {b}ound for {s}tatistical {l}earning.
\textit{Ann. Statist.} \textbf{34} 2326--2366.

\bibitem{n00}
Nemirovski, A. (2000).
Topics in non-parametric {s}tatistics.
\textit{ Ecole d'{\'E}t{\'e} de Probabilit{\'e}s de
Saint-Flour 1998. Lecture Notes in
Math}. \textbf{1738} 85--277.
New York: Springer.
\MR{1775640}

\bibitem{ssbook02}
Sch\"olkopf, B. and Smola, A. (2002).
\textit{Learning with Kernels}.
MIT Press.

\bibitem{ss05}
Steinwart, I. and Scovel, C. (2005).
Fast {r}ates for {s}upport {v}ector {m}achines.
In \textit{ Proceedings of the 18th Annual Conference on Learning Theory,
COLT 2005}. Berlin: Springer.
\MR{2203268}

\bibitem{ss04}
Steinwart, I. and Scovel, C. (2007).
Fast {r}ates for {s}upport {v}ector {m}achines using {G}aussian
{k}ernels.
\textit{Ann. Statist.} \textbf{35} 575--607.

\bibitem{tsy03}
Tsybakov, A.B. (2003).
Optimal rates of aggregation.
In \textit{Computational Learning Theory and Kernel Machines}
(B. Sch{\"o}lkopf and M. Warmuth, eds.). \textit{Lecture Notes in Artificial
Intelligence} \textbf{2777} 303--313.
Heidelberg: Springer.

\bibitem{tsy04}
Tsybakov, A.B. (2004).
Optimal aggregation of classifiers in statistical learning.
\textit{Ann. Statist.} \textbf{32} 135--166.
\MR{2051002}

\bibitem{v90}
Vovk, V.G. (1990).
Aggregating {s}trategies.
In \textit{Proceedings of the 3rd Annual Workshop on Computational
Learning Theory, COLT90}
371--386. San Mateo, CA: Morgan Kaufmann.

\bibitem{yang99}
Yang, Y. (1999).
Minimax {n}onparametric {c}lassification. I. {R}ates of
{c}onvergence.
\textit{IEEE Trans. on Inform. Theory} \textbf{45}
2271--2284.
\MR{1725115}

\bibitem{yang299}
Yang, Y. (1999).
Minimax {n}onparametric {c}lassification. II. {M}odel
{s}election for {a}daptation.
\textit{IEEE\break Trans. Inform. Theory} \textbf{45}
2285--2292.
\MR{1725116}

\bibitem{yang00}
Yang, Y. (2000).
Mixing strategies for density estimation.
\textit{Ann. Statist.} \textbf{28} 75--87.
\MR{1762904}

\bibitem{z04}
Zhang, T. (2004).
Statistical behavior and consistency of classification methods based
on convex risk minimization.
\textit{Ann. Statist.} \textbf{32} 56--85.
\MR{2051001}
\end{thebibliography}
\end{document}